\documentclass[a4paper,12pt,numbers,sort&compress]{article}

\pdfoutput=1

\usepackage[paper=letterpaper,margin=1in]{geometry}

\usepackage{amsmath,amsthm,amssymb,amsfonts,
epsfig,cite,setspace,bigstrut,longtable,array,breqn,url,color}

\usepackage{graphicx}
\usepackage[textwidth=3.4cm,textsize=10pt,disable]{todonotes}

\usepackage[singlelinecheck=false,justification=justified]{caption}

\usepackage{hyperref}
\usepackage{pdfpages,lipsum}

\renewcommand{\baselinestretch}{1.2}
\newcommand{\refpart}[1]{{\it (#1)}}  

\newcommand{\hpg}[5]{{}_{#1}\mbox{\rm F}_{\!#2}\!
  \left(\left.{#3 \atop #4}\right| #5 \right) }

\newcommand{\Pvi}{{\rm P}_{\rm VI}}

\begin{document}


\vskip 0.25in

\newcommand{\sref}[1]{\S~\ref{#1}}
\newcommand{\nn}{\nonumber}
\newcommand{\tr}{\mathop{\rm Tr}}
\newcommand{\comment}[1]{}
\newcommand{\cW}{{\cal W}}
\newcommand{\cN}{{\cal N}}
\newcommand{\cQ}{{\cal Q}}
\newcommand{\cH}{{\cal H}}
\newcommand{\cK}{{\cal K}}
\newcommand{\cZ}{{\cal Z}}
\newcommand{\cO}{{\cal O}}
\newcommand{\cA}{{\cal A}}
\newcommand{\cB}{{\cal B}}
\newcommand{\cC}{{\cal C}}
\newcommand{\cD}{{\cal D}}
\newcommand{\cE}{{\cal E}}
\newcommand{\cF}{{\cal F}}
\newcommand{\cG}{{\cal G}}
\newcommand{\cM}{{\cal M}}
\newcommand{\cX}{{\cal X}}
\newcommand{\IA}{\mathbb{A}}
\newcommand{\IP}{\mathbb{P}}
\newcommand{\IQ}{\mathbb{Q}}
\newcommand{\IH}{\mathbb{H}}
\newcommand{\IR}{\mathbb{R}}
\newcommand{\IC}{\mathbb{C}}
\newcommand{\IF}{\mathbb{F}}
\newcommand{\IV}{\mathbb{V}}
\newcommand{\II}{\mathbb{I}}
\newcommand{\IZ}{\mathbb{Z}}
\newcommand{\re}{{\rm~Re}}
\newcommand{\im}{{\rm~Im}}

\newcommand{\xx}{x} 

\newcommand{\ol}{\overline}
\newcommand{\diff}{\partial}
\newcommand{\dbar}{\ol{\partial}}

\newcommand{\tmat}[1]{{\tiny \left(\begin{matrix} #1 \end{matrix}\right)}}
\newcommand{\mat}[1]{\left(\begin{matrix} #1 \end{matrix}\right)}

\let\oldthebibliography=\thebibliography
\let\endoldthebibliography=\endthebibliography
\renewenvironment{thebibliography}[1]{%
\begin{oldthebibliography}{#1}%
\setlength{\parskip}{0ex}%
\setlength{\itemsep}{0ex}%
}%
{%
\end{oldthebibliography}%
}

\theoremstyle{theorem}
\newtheorem{theorem}{Theorem}[section]
\theoremstyle{definition}
\newtheorem{definition}{Definition}[section]
\theoremstyle{definition} 
\newtheorem{remark}[definition]{Remark} 
\theoremstyle{theorem}
\newtheorem{proposition}{Proposition}[section]

\def\thetheorem{\thesection.\arabic{proposition}}
\def\thetheorem{\thesection.\arabic{theorem}}
\def\thetheorem{\thesection.\arabic{definition}}

\setlength{\parskip}{0cm}
\setlength{\parindent}{0.5cm}
\setlength{\topsep}{0.35cm plus 0.2cm minus 0.2cm}

\def\theequation{\thesection.\arabic{equation}}
\newcommand{\setall}{\setcounter{equation}{0}
        \setcounter{theorem}{0}}
\newcommand{\setequation}{\setcounter{equation}{0}}
\renewcommand{\thefootnote}{\fnsymbol{footnote}}

~\\
\vskip 1cm

\begin{center}
{\Large \bf Genus One Belyi Maps by \\ Quadratic Correspondences}
\end{center}
\medskip

\vspace{.4cm}
\centerline{
{\large Raimundas Vidunas}$^1$ \ \&
{\large Yang-Hui He}$^2$
}
\vspace*{3.0ex}

\renewcommand{\baselinestretch}{0.5}
\begin{center}
{\it
  {\small
    \begin{tabular}{cl}
      ${}^{1}$ 
      & Graduate School of Information Science and Technology, \\
      & Osaka University, Osaka, Japan.\\
      ${}^{2}$
      & School of Physics, NanKai University, Tianjin, 300071, P.R.~China \& \\
      & Department of Mathematics, City, University of London, EC1V 0HB, UK \& \\
      & Merton College, University of Oxford, OX14JD, UK\\
      &\\
      & \qquad
      {\rm rvidunas@gmail.com, hey@maths.ox.ac.uk}
    \end{tabular}
  }
}
\end{center}

\renewcommand{\baselinestretch}{1.2}

\vspace*{4.0ex}
\centerline{\textbf{Abstract}} \bigskip
We present a method of obtaining a Belyi map on an elliptic curve from that on the Riemann sphere.
This is done by writing the former as a radical of the latter, which we call a quadratic correspondence, with the radical determining the elliptic curve.
With a host of examples of various degrees we demonstrate that the correspondence is an efficient way of obtaining genus one Belyi maps.
As applications, we find the Belyi maps for the dessins d'enfant which have arisen as brane-tilings in the physics community, including ones,
such as the so-called suspended pinched point, which have been a standing challenge for a number of years.

\newpage

\renewcommand{\baselinestretch}{0}
\tableofcontents
\renewcommand{\baselinestretch}{0.5}

\vspace{2cm}

\section{Introduction and Summary}\setall

There is a growing interest in, and therewith, an increasing demand for, explicit Belyi \cite {belyi} maps,
which are rational surjections onto $\IP^1$ from a Riemann surface $\Sigma$ ramified at exactly 3 points.
Such a demand arises from both mathematics, in the context of further understanding
dessins d'enfant \cite{catalog,leila,jones,lz,gg,jw},
Painlev\'e equations \cite{Kit1,KV1,LiTy}
and explicit computations \cite{hv,kmsv,sv,Vidunas:2016xun,vk},
and from physics, in the context of dessin realizations 
\cite{Jejjala:2010vb,Hanany:2011ra,Vidunas:2016xun} of brane-tilings in string theory \cite{Franco:2005rj,Franco:2005sm,Hanany:2005ve,Feng:2005gw,Feng:2000mi} (cf.~the catalogues of \cite{Davey:2009bp,Franco:2017jeo}), 
Bogomol'nyi-Prasad-Sommerfield (BPS) quivers
\cite{He:2015vua,Gabella:2017hpz} and Seiberg-Witten curves \cite{Ashok:2006br}.

A chief motivation comes from the fact that the bipartite graphs arising from the mathematics and the physics seem to reside in different worlds: the fundamental transformations such as Galois conjungation for the dessins, or Seiberg duality/cluster mutation for the quiver gauge theories associated to the tiling, etc., at first glance appear unnatural and mysterious to the other world. Yet, any translation between them would be of profound interest and in order to do, we need to collect data by having as many explicit Belyi maps as possible.

Computing Belyi maps of positive genus is a difficult problem \cite{sv}, for which
reduction to a computation of genus 0  Belyi maps helps enormously.
This is possible when the dessins have particular symmetries. 
This was the theme of \cite{Vidunas:2016xun}, 
where the list of genus 1 Belyi maps for known tilings was significantly expanded
by composing simple covering from elliptic curves to the Riemann sphere 
with simpler (univariate) genus zero Belyi maps, or with isogenies between elliptic curves.

In this paper, we continue with this theme of obtaining genus 1 Belyi maps 
from genus 0,  but introduce 
a different approach.
We shall use a quadratic transformation inspired by a common symmetric representation
of algebraic Painlev\'e VI functions, which we will call {\em quadratic correspondence}.
This simple transformation efficiently gives many more explicit Belyi maps.
In particular, we continue to address the catalogue 
of \cite{Davey:2009bp,Franco:2017jeo}, especially the former, which contains
many of the brane-tilings that have become standard to the literature.

We obtain the explicit Belyi maps for tilings \cite[(2.2),\,(3.4)--(3.6),\,(3.13)]{Davey:2009bp},
corresponding to the affine toric Calabi-Yau threefold geometries 
{\em SPP},  $L^{131}$, $L^{222}(II)$, $dP_1$, $L^{232}$
(see the table on p9 of \cite{He:2016fnb} for an explanation of the notation and context within the landscape of Calabi-Yau geometries).
These are recalled as Appendix tilings \refpart{B1}, \refpart{E1}, \refpart{C},  
\refpart{E2}, \refpart{G1} respectively in this paper.
This is assuring, especially given that {\em SPP} has been a standing challenge.
We thus cover all the applicable cases of \cite{Davey:2009bp} constructible by the quadratic correspondence.

We should also emphasize that there are many cases therein which are marked as ``inc.'', signifying that they are
inconsistent as physical brane-tilings (in that the corresponding quiver and superpotential indeed can be geometrically engineered as
a supersymmetric gauge theory whose moduli space is a Calabi-Yau 3-fold embeddable into string theory).
Nevertheless, these tilings, as purely abstract bipartite graphs on a torus, are well-defined and it is of undoubted
mathematical interest to find the explicit Belyi maps/dessins d'enfants for completeness; this we also do for a number of examples.

The paper is organized as follows.
We begin with the definition of the quadratic correspondence in \S\ref{sec:correspond}, explaining its
passport structure and in how to proceed from the Belyi map from genus 0 to genus 1.
Then, in \S\ref{sec:eg}, we provide a multitude of explicit examples, 
ranging from degree 7 to 12.
Finally, in \S\ref{sec:painleve}, we discuss the 
\todo{Reworded again.
"Commonality", "common representation" might be confusing, as the meanings differ}
commonality with the mentioned representation of algebraic Painlev\`e VI solutions, 
which originally inspired our quadratic correspondence.

\subsection*{Nomenclature}
\todo{Shifted a sentence here. The notations will be the same eventually.} 
Throughout we will continue with the notation of \cite{Vidunas:2016xun}.
\begin{itemize}
\item We use 
the multiplicative notation
  $[a_1^{p_1}\cdots a_{k_0}^{p_{k_0}} / b_1^{q_1}\cdots b_{k_1}^{q_{k_1}} /
  c_\infty^{r_1}\cdots c_{k_0}^{r_{k_\infty}}]$ for the branching passport
  of Belyi maps. This gives the branching
  indices, with repetition written as exponents, 
  of the pre-images of $0,1,\infty$ respectively.
\item By genus $g$ map we mean a Belyi map from a genus $g$ Riemann surface onto $\IP^1$.
\item A genus 1 curve is given as elliptic curve $\cE$ in a Weierstra\ss\ form 
$Y^2 = X^3+\ldots$ as {\em elliptic curves} 
$(\infty, \infty)$  as the origin of the group law.
The Klein J-invariant is denoted as $j$ and we adhere to the normalization that for $y^2 = x^3-x$, $j=1728$  (not $j=1$).
\item 
%
The elliptic curve $\cE$ is drawn as a parallelogram (the fundamental region)
with opposite sides identified, forming a doubly-periodic tiling of the (complex) plane.
The dessins d'enfant of genus 1 Belyi maps from $\cE$
are thus given as doubly periodic bipartite tilings, 
with black and white nodes representing the points above $\infty,0$
(respectively), and the cells represent the points above $1$.
\item If the $j$-invariant of $\cE$ is real, we draw 
its fundamental region in a doubly periodic
tiling as a rectangle (if the curve can have two real components)
or as a rhombus (if it can have one real component). 
If additionally the Belyi map is defined over $\IR$, the drawn tiling will have reflexion symmetries representing the complex conjugation:  
two reflexions parallel to a side of a rectangular fundamental region, 
or a reflexion parallel to a diagonal of a rhombus. 
In this way, our pictures reflect basic symmetries of the elliptic curves and the dessins.
In particular, a square fundamental domain signifies an elliptic curve with $j=1728$.

\comment{
\item We let $\omega_n=\exp(2\pi i/n)$ denote the primitive $n$-th root of unity.

\todo{I think we need thisffollowing def because we refer to it late in first paragraph of Sec 2.1}
\item The Gau\ss{} hypergeometric function is defined via the Euler integral: 
  \[
  \hpg{2}{1}{a,b}{c}{z} = \frac{\Gamma(c)}{\Gamma(b)\Gamma(c-b)}
  \int_0^1 t^{b-1}(1-t)^{c-b-1}(1 - tz)^{-a} dt.
  \]
It satisfies the hypergeometric differential equation 
\begin{equation} \label{eq:HGE}
z\,(1-z)\,\frac{d^2y(z)}{dz^2}+
\left(c-(a+b+1)\,z\right)\,\frac{dy(z)}{dz}-a\,b\,y(z)=0.
\end{equation}  
This is a canonical Fuchsian equation with 3 singularities, namely $z=0$, $z=1$, $z=\infty$.
Alternative symmetric parameters are the {\em local exponent differences} $1-c,c-a-b,a-b$
at the three singular points.

}

\end{itemize}

\section{Correspondence Reduction to Genus Zero}\setall
\label{sec:correspond}

The basic construction of this paper is based on the following
operation on a rational or algebraic function $\varphi_0$:
\begin{align} \label{eq:basicop}
\varphi_1 =
\cQ\{\varphi_0\}:=&\, \frac12+\frac{\sqrt{1-\varphi_0}}2.
\end{align}
Geometrically, we have a 2-to-2 algebraic correspondence 
between the coverings $\varphi_0$ and $\cQ\{\varphi_0\}$.
The genus of $\cQ\{\varphi_0\}$ may be higher than that of $\varphi_0$.
\todo{Reworded again, as quadratic Painleve stuff is gone.}
This type of radical expression has been used in representing algebraic
Painlev\'e VI solutions, as we 
\todo{ discuss $\implies$ recall, recollect?} 
discuss in \S \ref{sec:painleve}.


We are particularly interested in the cases when $\varphi_0$ and $\varphi_1$ 
are Belyi functions of genus 0 (as a Riemann sphere $\IP^1$) and 1 (as an elliptic curve $\cE$).
The branching passport of  $\varphi_1$ 
will have the same branching pattern in the fibers $\varphi_1=0$ and $\varphi_1=1$, where
$\sqrt{1 - \varphi_0} = \mp 1$.
Correspondingly, the dessins d'enfant of these  $\varphi_1$ 
will have an involution symmetry between black and white points.
This is the characteristic property of genus 1 Belyi maps that can be 
reduced to genus 0 Belyi maps by the $\cQ$-correspondence.
The involution symmetry gives a degree 2 covering $\cE\to\IP^1$.
The presumed elliptic curve $\cE$ is defined by the square root $\sqrt{1-\varphi_0}$.

\paragraph{Passport Structure: }
Given $\varphi_1=\cQ\{\varphi_0\}$, possible passports for 
$\varphi_0=4\varphi_1(1-\varphi_1)$ can be easily determined.
The applicable passport of $\varphi_1$ has the schematic form $[A/A/B^2C]$,
where $A,B,C$ are some partitions of integers.
The anticipated passport of $\varphi_0$ is $[A/2^{\ell}1^k/\widetilde{B}\,C]$
with $k\le 4$. 
This enforces the following restrictions 
\begin{enumerate}
\item 
The partition $C$ contains $4-k$ numbers, 
and the partition $\widetilde{B}$ has every number of $B$ doubled.
\item The covering $\cE\to\IP^1$
branches over the points represented by $1^k$ and $C$. 
\item The branching orders in $A$ represent 
the black points above $\varphi_0=0$,
and the pairs of anti-symmetric points in $\varphi_1=0$, $\varphi_1=1$.
\item With the exception of $k$ non-branching points,
the white points in $\varphi_0=1$ have the branching order 2.
\item $\widetilde{B}$ represents the $\varphi_0$-cells
where $\cE\to\IP^1$ does not branch, 
while $B$ represents the corresponding pairs of the $\varphi_1$-cells.
The valencies in $B$ and $\widetilde{B}$ differ by the factor 2.
\item In addition, we require that all valencies in $C$ be odd numbers.
\end{enumerate}
The last condition excludes the case when
the valencies in $1^k$ and $C$ are of mixed parity.
The desired relation between $\varphi_0$, $\varphi_1$ 
is not there then, as will be exemplified in Remark \ref{rm:parity2}. 
In particular, the square root $\sqrt{1-\varphi_0}$ defines
an algebraic curve of genus 0, and it is not a regular function 
on any genus 1 curve.

\begin{definition}
Given a genus 1 passport of the shape $[A/A/B^2C]$,
the genus 0 passports $[A/2^{\ell}1^k/\widetilde{B}\,C]$ described above
are called {\em correspondent passports}. 
\end{definition}
\begin{remark}
  \label{rm:parity1}
If $k=0$ and all valencies in $C$ are even, then 
$\cQ\{\varphi_0\}$ is minimally defined as a genus 0 function.
Genus 1 functions are then defined by the straightforward 
composition construction of \cite{Vidunas:2016xun}.
For example, the genus 0 map 
\begin{equation}  \label{eq:isog2o0}
\varphi_0=-\frac{(x^2-1)^4}{16x^2(x^2+1)^2}
\end{equation}
with the passport $[4^2/2^4/2^4]$ gives the Belyi map
\begin{align} \label{eq:isog2o}
\varphi_1  
=\frac{(x+1)^4}{8\,(x^3+x)} 
\qquad \mbox{as} \qquad \varphi_1  =\cQ\left\{\varphi_0\right\}
\end{align}
on $\cE$ being $y^2=x^3+x$, with the passport $[4^2/2^4/4^2]$.
Since $\cQ\left\{\varphi_0\right\}$ has two values, 
we can interpret $\cQ\left\{\varphi_0\right\}=1-\varphi_1$ as well.
As a rational function, $\varphi_1\in\IC(x)$ has degree 4.
The genus 1 map is related to a 2-isogeny on 
$y^2=x^3+x$ and the Belyi map $\psi_1=(x^2+1)^2/(4x^2)$ 
in \cite[(4.9)]{Vidunas:2016xun}, namely
\begin{equation}
\varphi_1=\frac{\psi_1}{\psi_1-1}=\frac{(x^2+1)^2}{(x^2-1)^2}
\end{equation}
after the substitution
\begin{equation}
x\mapsto \frac{(x+1)^2+2\sqrt{2} y}{(x-1)^2}.
\end{equation}

In a general situation when $\cQ\{\varphi\}$ is minimally a genus 0 function,
the genus 1 passport of $\cQ\{\varphi\}$ 
is not necessarily the targeted $[A/A/B^2C]$, as will be demonstrated in Remark \ref{rm:deg12g00}.
\end{remark}

\subsection{Computational Methods}

The genus 0 maps are obtained by using the {\sf Maple} program developed by the first author in collaboration with Mark van Hoeij in \cite{vHVHeun}. 
The program can effectively compute Belyi maps of genus 0 if the passport is nearly regular.
Given positive integers $k,\ell,m,n$,
a Belyi map in \cite{vHVHeun} is defined to be {\em $(k,\ell,m)$-minus-$n$ regular} 
if, with exactly $n$ exceptions in total, all points above $z=1$ have branching order $k$, 
all points above $z=0$ have branching order $\ell$, 
and all points above $z=\infty$ have branching order $m$. 
Typically, these Belyi maps can pull-back the differential equation for the 
hypergeometric function
$\hpg21{\frac12(1-\frac1k-\frac1\ell-\frac1m),\,\frac12(1-\frac1k-\frac1\ell+\frac1m)}
{1-1/k}{z}$
to Fuchsian equations with $n$ singularities (cf.~\cite{vf,vHVHeun}). 
When $n$ is small, this gives additional algebraic equations 
for the undetermined coefficients of the desired Belyi maps.
See \cite[\S 5.2]{vHVHeun}, \cite{hv} for the full system of utilized relations
of undertermined coefficients.

In \cite{vHVHeun}, 366 Galois orbits of Belyi $(k,\ell,m)$-minus-$4$ 
Belyi maps of degree $\le 60$
were computed. 
Not all computed Belyi maps of genus 0 in the present article are highly regular, 
as we have $n\le 6$. But their degree is at most 12, 
so they can be computed within minutes.


Having defined the type of Belyi maps on which we focus, it is expedient to point
out some of details in our computation, especially with regard to the covering of
the elliptic curve $\cE$ over $\IP^1$.
Now, simplification of computed Belyi maps (by M\"obius transformations
or isomorphisms of their genus 1 curve) is a cumbersome problem,
which is still not satisfactorily automate in general. 
In particular, bringing an equation for a genus 1 curve to 
the most canonical Weierstra\ss\ form $y^2 = 4x^3 - g_2 x - g_3$ 
is not always practical.
A transformation from the genus 1 form
$y^2=x^4+ax^3+bx^2+cx+d$ is presented in \cite[\S 2.1]{Vidunas:2016xun}.

\comment{
Often, our construction from a genus 0 map 
directly gives elliptic curves in the form
\todo{For these formulas, we can refer to our first paper.}
\begin{align}
y^2=x^4+ax^3+bx^2+cx+d 
\ .
\end{align}
By keeping a point at infinity, we can transform this curve to the Weierstra\ss\ form 
\begin{equation} \label{eq:ecwfrom4}
Y^2=X^3+bX^2+(ac-4d)X+(a^2-4b)d+c^2
\end{equation}
by the transformation
\begin{align} \label{eq:ec4to3}
x &= -\frac{2Y+aX+2c}{4X-a^2+4b}, \\
\nn
y &=  \frac{8cY-4X^3+4(ac-4d)X+8c^2+(a^2-4b)(aY+3X^2+2bX+a c+4d)}{(4X-a^2+4b)^2}
\ .
\end{align}
}

\begin{remark}
  {
Non-existence of genus 0 Belyi maps with a given passport can be easily proven 
by contradiction to implied pull-back transformations to non-existent
Fuchsian equations \cite[\S 5]{vf}, \cite[Remark 4.1]{vHVHeun}. 
In particular, this excludes many $(k,\ell,m)$-minus-$n$ regular maps with $n\le 2$.
Among our examples, there are no Belyi maps with the genus 0 passports
$[5\,3/2^4/2^4]$, 
$[5\,4\,3/2^6/2^4\,4]$. 
However, a $(3,2,2)$-minus-2-regular 
Belyi map with the passport $[3^2\,6/2^4\,4/2^6]$ does exist, 
as presented in Remark \ref{rm:deg12g00}. 
It gives a pull-back transformation of the differential equation 
for $\hpg21{\!-1/6,\,1/3}{2/3}{z}$ 
to a Fuchsian equation with two apparent singularities.
  }
\end{remark}

\section{Explicit Examples}\setall
\label{sec:eg}
Using the general methods outline above, we now demonstrate with a multitude of examples, of various degree.
We construct genus 1 Belyi maps with the passports 
\begin{align}
\begin{matrix} 
[4\,3/4\,3/2^2\,3], & [5\,3/5\,3/2^4], &  [3^3/3^3/2^2\,5], 
&  [3^24/3^24/3^22^2], \\
\; [4^2\,3/4^2\,3/2^4\,3], \; & \; [3^2\,5/3^2\,5/2^4\,3], \; & \;
[3^2\,6/3^2\,6/2^6], \; & \; [5\,4\,3/5\,4\,3/2^6]. \;
\end{matrix}
\end{align}
Their dessins d'enfant give the tilings 
\cite[(2.2),\,(2.3),\,(2.6),\,(3.3)--(3.7),\,(3.9),\,(3.10), 
(3.13),\;(3.16),\;(3.18)\,--\,(3.21),\;(3.24),\;(3.29)\,--\,(3.33),\;(3.35)]
{Davey:2009bp}. 
The theories 
\cite[(2.1),
\,(2.4),\,(2.5),\,(3.1),\,(3.2),\,(3.26),\,(3.27)]{Davey:2009bp} are so-called ``orbifold'' which 
we leave to an upcoming paper to address systematically in relation to isogenies.
There are a few other tilings in \cite{Davey:2009bp} that have the considered
or similarly symmetric passports. 
Inapplicability of the correspondence construction to them
is explained in Remark \ref{rm:except}.

\subsection{A Belyi map of degree 9}
\label{sec:d9}

We start with this genus 0 map of degree 9, with the passport $[3^3/2^3\,1^3/4\,5]$:
\begin{align}
\psi_9= &  -\frac{\big(x^3+45x-450)^3}{2916\,\big(3x-5\big)^4} \\
= & \, 1-\frac{ \big(x^3-3x^2-9x+315\big)^2\,\big(x^3+6x^2+180x-900\big)}
{2916\,\big(3x-5\big)^4}.  \quad \nonumber
\end{align}
We have
\begin{equation} \label{eq:f26}
\cQ\{\psi_9\} =\frac12+\frac{x^3-3x^2-9x+315}{108\,\big(3x-5\big)^2}
\,\sqrt{x^3+6x^2+180x-900}\,.
\end{equation}
Let $\cE_9$ denote the elliptic curve $y^2=x^3+6x^2+180x-900$.
The function $\cQ\{\psi_9\}$  is a genus 1 map on $\cE_9$, 
of degree 9 as well, with the passport $[3^3/3^3/2^2\,5]$. 
The $j$-invariant of $\cE_9$ 
equals $6\,(112/25)^3$.

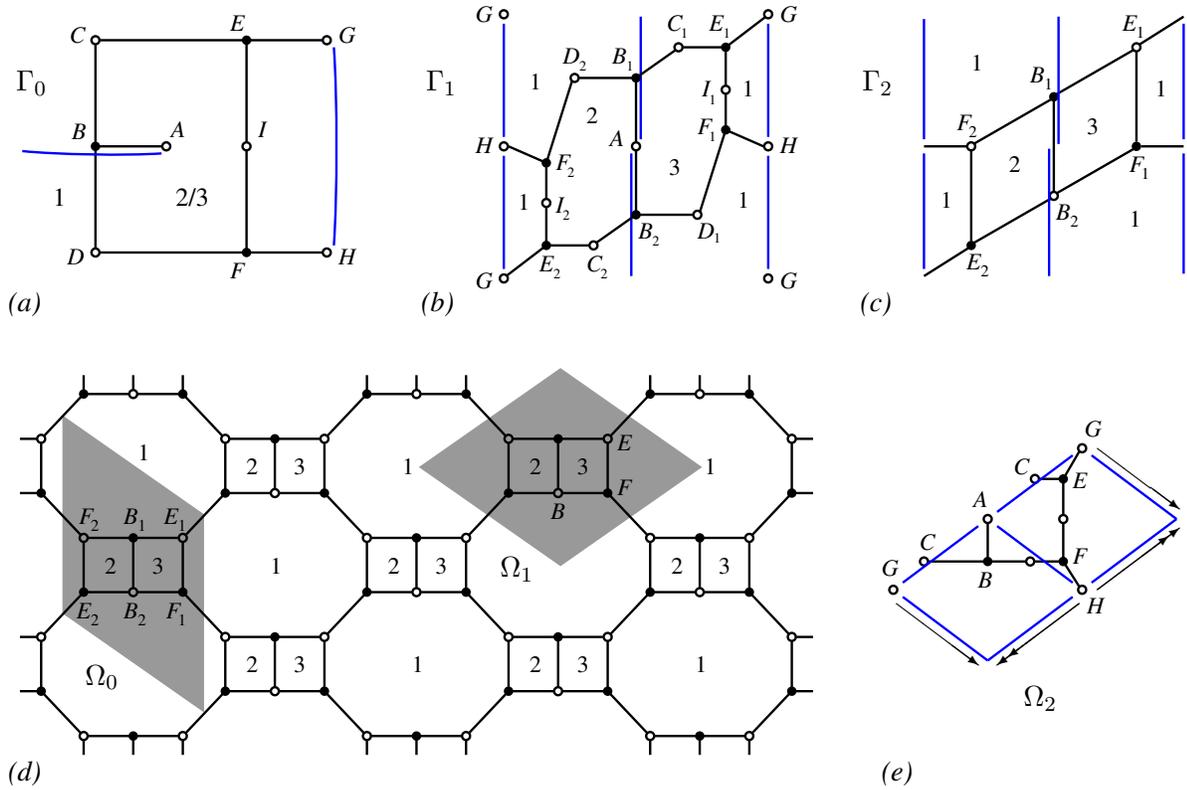
\begin{figure}[t] 
\begin{picture}(460,300)
\put(-2,6){\includegraphics[width=450pt]{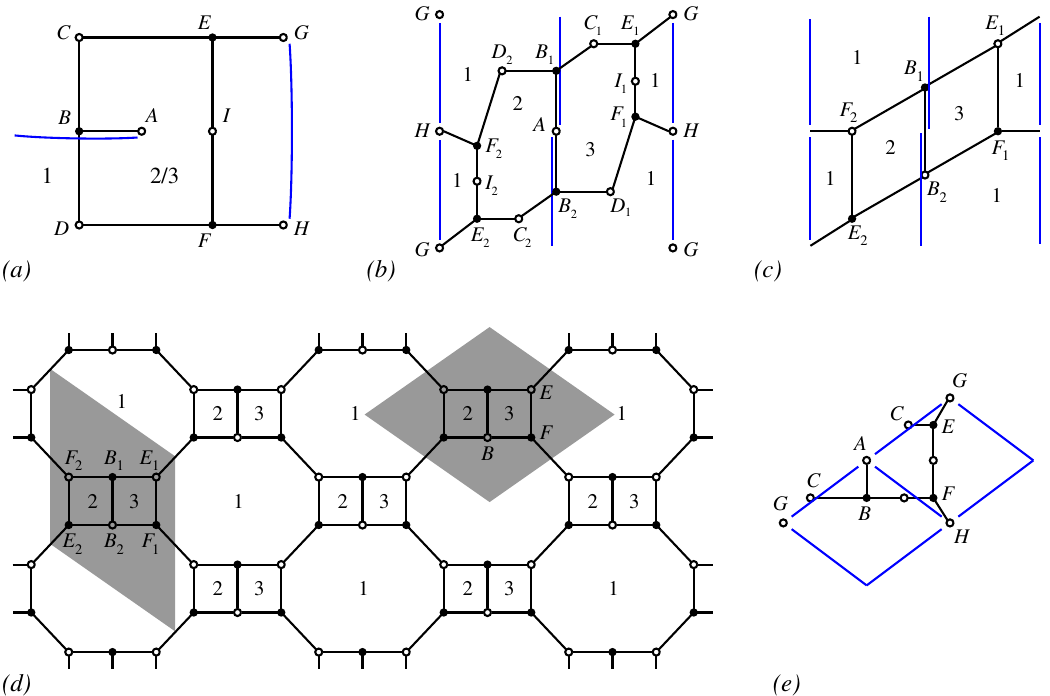}}
\put(2,269){\small $\Gamma_0$}
\put(157,269){\small $\Gamma_1$}
\put(322,269){\small $\Gamma_2$}
\put(28,44){\small $\Omega_0$}
\put(185,85){\small $\Omega_1$}
\put(383,36){\small $\Omega_2$}
\put(335,75){\vector(4,-3){31}}  \put(411,133.75){\vector(4,-3){31}}
\put(411,80.25){\vector(4,3){31}}  \put(404,75){\vector(-4,-3){31}}
\put(411,80.25){\vector(4,3){27}}  \put(404,75){\vector(-4,-3){27}}
\end{picture}
\caption{From a genus 0 map to a corresponding tiling \cite[(3.3)]{Davey:2009bp}}
\label{fig:dessin33}
\end{figure}

This  {\em dessin d'enfant} of $\cQ\{\psi_9\}$ is equivalent 
to the tiling \cite[(3.3)]{Davey:2009bp}.
It can be obtained in the following steps:
\begin{enumerate}
\item Draw the genus 0 dessin $\Gamma_0$ 
of $\psi_9$ as in Figure \ref{fig:dessin33}\refpart{a}.
The branching points of the projection $x:\cE_9\to\IP^1$ 
are going to be $A,G,H$ and $\infty$. We mark the cuts $GH$ and $A\infty$.
The inner cell will have two pre-images ``2" and ``3" on $\cE_1$, 
thus we mark it by ``2/3".  The outer cell will have one pre-image ``1", 
because it contains the branching point $\infty$. 

\item Draw the dessin $\Gamma_1$ of $\psi_9$ as a Belyi function on $\cE_9$. 
The elliptic curve $\cE_9$, as always, is represented by a fundamental domain of a 
doubly periodic lattice in $\IC$. 
Since the dessin is defined up to homotopy, we draw 
the fundamental domain as a rectangle initially. In Figure \ref{fig:dessin33}\refpart{b},
the point $G$ is represented by the corners of the rectangle.
The other 3 branching points are placed regularly so that the cuts $GH$, $A\infty$ 
and the complementary cuts $AH$, $G\infty$ give a subdivision into 4 smaller rectangles 
--- call them {\em quads}.
Each of the quads represents (alternatively) either the inside or the outside 
of a curvilinear rectangle {\em GHA$\infty$} around $\Gamma_0$ 
in Figure \ref{fig:dessin33}\refpart{a}. 
The dessin $\Gamma_1$ is drawn as a double cover of $\Gamma_0$,
so that a sketch in each quad traces 
the corresponding half {\em GHA$\infty$} of $\Gamma_0$. 
Figure \ref{fig:dessin33}\refpart{b} gives the corresponding labels of vertices and cells,
and keeps the central symmetry of $\Gamma_1$.
This is the same routine as outlined in \cite[\S 2.2]{Vidunas:2016xun}.
The central symmetry is a characteristic feature of dessins d'enfant
(within a proper fundamental domain) of the composite Belyi maps in
\cite{Vidunas:2016xun}.

\item All white vertices of $\Gamma_1$ have the valency 2.
To modify $\Gamma_1$ to the dessin $\Gamma_2$ of $\cQ\{\psi_9\}$,
the white vertices are discarded, then the black vertices are newly bi-coloured 
black and white so that a bipartite graph is obtained.  Figure \ref{fig:dessin33}\refpart{c}
depicts a straightened-up version of $\cQ\{\psi_9\}$.
The dessin has a central {\em anti-symetry} between the white and black points.
\item Copies of $\Gamma_2$ in a fundamental domain 
can be stacked to a doubly periodic tiling  of the complex plane. 
Since $\cE_9$ has one $\IR$-component, an appropriate shape
of the fundamental domain is a rhombus. (It would remain a rectangle 
for elliptic curves  with two $\IR$-components.) 
Figure \ref{fig:dessin33}\refpart{d} depicts $\Gamma_2$
as a doubly-periodic tiling. The straightforward picture is obtained 
by stretching the fundamental domain to the rhombus $\Omega_0$. 
The other fundamental domain $\Omega_1$ exhibits the 
$\IR$-involution symmetry (by the vertical reflection), 
and the horizontal and central anti-symmetries between the black and white vertices.
\end{enumerate}
\begin{remark} 
  \label{rm:noncorr}
A passport symmetry between branching orders of the black and white points 
of a genus 1 dessin is not sufficient for the correspondence reduction to a genus 0 map.
Appendix figure \refpart{M} 
gives two counter-examples 
that we mention again in Remark \ref{rm:except} and \S \ref{sec:notp6}.
A defining criterium for the correspondence reduction
is existence (up to homotopy) of a central anti-symmetry 
between the black and white points.

Given a genus 1 dessin with a central anti-symmetry,
a correspondent genus 0 dessin is obtained as follows:
\begin{itemize}
\item Choose a fundamental domain with a central anti-symetry,
such as $\Omega_1$ in Figure \ref{fig:dessin33}\refpart{d}.
\item Take a proper half $\Omega_2$ of the fundamental domain 
(formed by two adjacent quads), 
and glue its boundary $\partial\Omega_2$ to a Riemann sphere
following the arrows and vertex labels in Figure \ref{fig:dessin33}\refpart{e}. 
\item Re-colour all vertices to black, and add white vertices to:
\begin{itemize}
\item the loose edge points on the glued boundary $\partial\Omega_2$;
\item the middle points of the edges 
between black points within  $\Omega_2$.
\end{itemize}
\end{itemize}
\end{remark}

\subsection{Maps of degree 8}
\label{sec:d8}

Next, we look for genus 1 degree 8 Belyi maps with the passport $[5\,3/5\,3/2^4]$.
There are two correspondent genus 0 passports: $[5\,3/2^21^4/4^2]$ and $[5\,3/2^4/2^4]$.
There are no genus 0 Belyi maps with the latter passport. 
That would have given an alternative example to Remark \ref{rm:parity1}.

There is one Belyi map 
with  the correspondent passport $[5\,3/2^21^4/4^2]$:
\begin{align}
\psi_8= \frac{64\,\big(4x-5)^3}{\big(x^2-10\big)^4}
= 1-\frac{ \big(x^2-2x+6\big)^2\,\big(x^4+4x^3-40x^2-200x+500\big)}
{\big(x^2-10\big)^4}. 
\end{align}
It leads to the genus 1 map
\begin{equation} \label{eq:f33}
\cQ\{\psi_8\} =\frac12+\frac{x^2-2x+6}{2\,\big(x^2-10\big)^2}
\,\sqrt{x^4+4x^3-40x^2-200x+500}.
\end{equation}
The genus 1 curve is defined by the square root.
The $j$-invariant equals $-5000$, as in \cite[(3.6)]{Vidunas:2016xun}. 
An isomorphic elliptic curve the Weierstrass form is
$Y^2=X^3+X^2-8X+8$. 
Appendix figure \refpart{A} depicts the genus 0 dessin of $\psi_2$ (with branch cuts marked)  
and the genus 1 dessin of $\cQ\{\psi_2\}$.
These dessins are related by the routine \refpart{i}--\refpart{iv} of the previous section.
The dessin d'enfant for \eqref{eq:f26} can be identified as  \cite[(2.6)]{Davey:2009bp}.

\begin{remark} 
  \label{rm:parity2}
It may seem that Belyi maps with the passport $[5\,3/5\,3/2^4]$ can be obtained
from genus 0 maps with the passport $[5\,3/2^31^2/2^24]$,
by allowing $e:\cE\to\IP^1$ to branch above the $2^2$ and $1^2$ points in
the correspondence construction. 
There is one Galois orbit of these genus 0 Belyi maps, defined over $\IQ(\sqrt{10})$:
\begin{equation}
\psi_3=\frac{(-25-\sqrt{10})(8x+3-3\sqrt{10})^3}{125x^4(x^2+4x-2-2\sqrt{10})^2}.
\end{equation}
The anticipated elliptic curve is defined by $y^2 = P_1P_2$ with 
\begin{equation}
P_1=x^2+4x-2-2\sqrt{10}, \qquad
P_2={\textstyle x^2+\Big(2+\frac{6\sqrt{10}}5\Big)x-1-\frac{4\sqrt{10}}{5}}.
\end{equation}
However, $\cQ\{\psi_3\}$ would involve $\sqrt{P_2}$ but not $\sqrt{P_1}$ or $\sqrt{P_1P_2}$.
Therefore $\cQ\{\psi_3\}$ would not be a regular function on $\cE$.
Appendix figure \refpart{N} 
demonstrates that we can draw the dessin 
of $\psi_3$ on $\cE$ (with the white vertices {\em cleanly} discarded), 
but the bi-colouring step \refpart{iii} for the anticipated dessin 
of $\cQ\{\psi_3\}$ is impossible. 
For example, equivalent points in adjacent fundamental domains are connected
by a sequence of adjoining edges of odd length 3.
\end{remark}

\subsection{Maps of degree 7}
\label{sec:d7}

Genus 1 maps with the passport $[4\,3/4\,3/2^23]$
can be obtained from genus 0 map with the passport $[4\,3/2^21^3/4\,3]$.
There is one Galois orbit of these maps, defined over $\IQ(\sqrt7)$, with the genus 0 being
\begin{align}
\psi_4= & \, -\frac{(9x-4\sqrt7-29)^3\,(x-2\sqrt7-4)^4}{(4\sqrt7+29)\,(7x-10\sqrt7-32)^4} \ ,
\end{align}
and the genus 1 map, 
\begin{equation} \label{eq:f22}
\cQ\{\psi_4\}=\frac12+\frac{27x^2-(72\sqrt7+252)x+344\sqrt7+1036}
{(4\sqrt7+2)\,\big(7x-10\sqrt7-32)^2}\,y
\end{equation}
on the curve
\begin{equation}
y^2=x\,\big(x^2-(4\sqrt7+7)x+16\sqrt7+35\big).
\end{equation} 
The associated dessin is identified with  \cite[(2.2)]{Davey:2009bp}. 
The genus 0 and 1 dessins are depicted in Appendix figures \refpart{B1}--\refpart{B2}.

\subsection{Maps of degree 10}

Here we look for genus 1 Belyi maps with the passport $[3^24/3^24/3^22^2]$.
There are two correspondent passports, 
namely $[3^24/2^41^2/3^24]$ and $[3^24/2^31^4/6\,4]$.
Let us present these in turn.

The genus 0 passport $[3^24/2^41^2/3^24]$ has two Galois orbits of Belyi maps, 
defined over $\IQ$ or $\IQ(\sqrt{-5})$.
The $\IQ$-rational 
map is
\begin{align}
\psi_{10}= &\, \frac{(9x-1)^3\,(x-9)^3\,(x+1)^4}{800000\,x^3\,(x-1)^4}.
\end{align}
The genus 1 map
\begin{align} \label{eq:comb33a}
\cQ\{\psi_{10}\} =
\frac12+\frac{27x^4-72x^3+602x^2-72x+27}{4000\,x^2\,(x-1)^2}
 \ \sqrt{-5x(x^2-18x+1)}
\end{align}
is equivalent to the map $\varphi_5$  in \cite[(2.13)]{Vidunas:2016xun}
on the curve on $y^2\!=\!x(x^2+18x+1)$.
One can transform 
this new expression to $\varphi_5/(\varphi_5-1)$
by $x\mapsto \big(\!-x^2+4\sqrt{5}\,y-38x-1\big)/(x-1)^2$.
\todo{Added this.}
The genus 1 dessin is identified with the brane tiling \cite[(3.5)]{Davey:2009bp}.
As illustrated by Appendix Figure \refpart{C} and \cite[Appendix \refpart{B}]{Vidunas:2016xun},
it is obtainable from genus 0 dessins by both correspondence and composition methods.
The geometric reason is that the tiling has both central anti-symmetries
(say, centered at the cells 1, 4) and central symmetries (say, centered at the cells 2, 4).
The degree of the genus 0 dessin in \cite{Vidunas:2016xun} 
equals 5 rather than 10, and the re-colouring step \refpart{iii} of \S\ref{sec:d9} 
does not apply there.



The $\IQ(\sqrt{-5})$-orbit gives the genus 1 Belyi map
\begin{align} \label{eq:qm310}
\cQ\left\{
\frac{(4x^2-\big(2+10\sqrt{-5})x-31+5\sqrt{-5}\big)^3\,(x+6)^4}
{50\,(5-\sqrt{-5})\,\big(x^2+5+2\sqrt{-5}\big)^3\,\big(x-2\sqrt{-5}\big)^4}
\right\}
\end{align}
on the curve
\begin{equation}
w^2=\big(15+12\sqrt{-5}\big)\,\big(x^2+5+2\sqrt{-5}\big)\,
\big(x^2+2x+16-2\sqrt{-5}\big) \ .
\end{equation} 
This curve is isomorphic to $y^2=x\,(x-1)\,\big(x+4\sqrt{-5}\big)$ over $\IC$.
Representative genus 0 and 1 dessins are given in Appendix \refpart{D}.
The complex conjugate tiling is \cite[(3.10)]{Davey:2009bp}.

The genus 0 Belyi maps with the other correspondent passport $[3^24/2^31^4/6\,4]$
are defined over the degree 5 number field 
\begin{equation} 
\IQ(\xi)/(\xi^5+10\xi^4+10\xi^3-20\xi^2-20\xi+16).
\end{equation}
Here is an expression: 
\begin{equation} \label{eq:deg10d5}
\psi^\star_{10}=\frac{4\,\big((\xi^3-2\xi+1)\,x^2-5\big)^{3}\,
\big(\frac12\xi^2\,x-5\big)^{4}}
{C\,\big(5x+\frac{25}4\xi^4+66\xi^3+99\xi^2-75\xi-180\big)^{4}}
\end{equation}
with $C=-\frac{18469}4\xi^4-\frac{12647}2\xi^3+\frac{18703}2\xi^2
+11459\xi-8443$. 
The genus 1 curve for $\cQ\{\psi^\star_{10}\}$ has a lengthy expression
$y^2=x^4+\ldots.$ It is isomorphic to
\begin{align}
y^2= &\; x^3+5\,(-231\xi^4-308\xi^3+474\xi^2+540\xi-407)\,x  \\
& +5\,(601119\xi^4+840844\xi^3-1220178\xi^2-1527900\xi+1118282).  
\quad \nonumber
\end{align}
The dessins are given in  Appendix figures \refpart{E1}--\refpart{E4}.
They are identified with \cite[(3.4),\,(3.6),\,(3.7),\,(3.9)]{Davey:2009bp},
respectively.

\subsection{Maps of degree 11}
\label{sec:deg11}

Here we consider genus 1 Belyi maps with these two passports:
$[4^2\,3/4^2\,3/2^4\,3]$ and $[3^2\,5/3^2\,5/2^4\,3]$. 
 The correspondent passports are, respectively,
 $[4^2\,3/2^4\,1^3/4^2\,3]$
 and $[3^2\,5/2^4\,1^3/4^2\,3]$. We again address these in turn.


There are two Galois orbits with the genus 0 passport  $[4^2\,3/2^41^3/4^2\,3]$.
One is defined over $\IQ(\sqrt{11})$:
\begin{equation} \label{eq:deg11r11}
\psi_{11}=\frac{67-44\sqrt{11}}{67+44\sqrt{11}} \,
\frac{\big(7x^2+\frac{11}{49}(81+2\sqrt{11})x+11\big)^4\,x^3}
{\big(11x^2+\frac{11}{49}(81-2\sqrt{11})x+7\big)^4}.
\end{equation}
The genus 1 function $\cQ\{\psi_{11}\}$
is defined on an elliptic curve $y^2=x^3+\ldots$ isomorphic to 
\begin{equation}
y^2=x^3+(24\sqrt{11}-88)\,x-172\sqrt{11}+561.
\end{equation}
The other Galois orbit is defined over the cubic field
$\IQ(\mu)/(\mu^3-\mu^2+4\mu+2)$. The genus 0 map is
\begin{equation}  \label{eq:deg11cubic}
\psi^0_{11}=\frac{\big(27x^2+11(\mu^2+5)x+11 (1-4\mu)\big)^4\,x^3}
{\big(11(1-4\mu)x^2+11(\mu^2+5)x+27\big)^4} \ .
\end{equation}
The genus 1 curve for $\cQ\{\psi^0_{11}\}$ is isomorphic to 
\begin{equation}
y^2=x\,(x^2+(8\mu^2-16\mu+35)\,x+88\mu^2-120\mu+403).
\end{equation}
There are several tilings in \cite{Davey:2009bp} 
with the passport $[4^2\,3/4^2\,3/2^4\,3]$.
To assign them to $\cQ\{\psi_{11}\}$ or $\cQ\{\psi^0_{11}\}$, 
we observe that the variable change $x\mapsto 1/x$ 
transforms $\psi^0_{11}$ to $1/\psi^0_{11}$,
but it transforms $\psi_{11}$ to the $\sqrt{11}$-conjugated $1/\psi_{11}$.
Hence the genus 0 dessins of  $\psi^0_{11}$ have the duality symmetry
between the cells and the black vertices, while this duality interchanges
the two dessins of $\psi_{11}$.
The genus 0 and 1 dessins defined over $\IQ(\sqrt{11})$
are depicted in Appendix figures \refpart{F1}--\refpart{F2}.
The tilings are identified as \cite[(3.18), (3.19)]{Davey:2009bp}.
The dessins defined over the cubic field
are depicted in Appendix figures \refpart{G1}--\refpart{G2}, up to complex conjugation.
The tilings are identified as \cite[(3.13), (3.20)]{Davey:2009bp}.


The second genus 0 passport $[3^2\,5/2^4\,1^3/4^2\,3]$
gives a Galois orbit defined over the degree 5 number field 
\begin{equation}
\IQ(\lambda)/(\lambda^5+3\lambda^4+4\lambda^3-10\lambda^2+26\lambda-54). 
\end{equation}
Here is an expression:
\begin{equation} \label{eq:deg11d5}
\psi^\star_{11}=
\frac{32\,(914\lambda^4+572\lambda^3-7010\lambda^2-25817\lambda+51759)\,H^3}
{115625\,x^3\, \big(x^2-11x+\frac{11}{370}
(14\lambda^4+49\lambda^3+25\lambda^2-331\lambda+1031)\big)^4},
\end{equation}
where
\begin{equation*} \textstyle
H=11x^2-\frac{11}{296} (14\lambda^4+86\lambda^3+136\lambda^2-257\lambda+1771)x+\frac12(\lambda^4+11\lambda^3+12\lambda^2-66\lambda+112).
\end{equation*}
The map $\cQ\{\psi^\star_{11}\}$ has the passport $[3^2\,5/3^2\,5/2^4\,3]$.
Its genus 1 curve 
is isomorphic to
\begin{align}
y^2= &\; 37x^3-3\,
(15758\lambda^4+70508\lambda^3+169511\lambda^2+99392\lambda+553828)\,x \\
& +2\,(1769584\lambda^4+7967657\lambda^3+19126222\lambda^2
+11359204\lambda+62899120).  \quad \nonumber
\end{align}
There is one real dessin and two pairs of complex conjugate dessins.
They are depicted in Appendix figures \refpart{H1}--\refpart{H3}.
They are identified with \cite[(3.16),\,(3.21),\,(3.24)]{Davey:2009bp},
respectively.

\begin{remark} 
  \label{rm:except}
The tiling \cite[(3.25)]{Davey:2009bp}
has the considered passport $[3^25/3^25/2^43]$, 
but it cannot be obtained by the correspondence construction.
In agreement with Remark \ref{rm:noncorr},
the tiling has no central anti-symmetry between the black and white points;
see the left tiling in Appendix figure 
\refpart{M}.

Similarly, the degree 12 tiling \cite[(3.37)]{Davey:2009bp} 
has the passport $[3^26/3^26/2^6]$ considered in \S \ref{sec:deg12b}.
The dessin has no central anti-symmetry, but there are central symmetries.
Hence it can be obtained by the composition construction 
in \cite[Fig.~3\refpart{b}]{Vidunas:2016xun}.

The degree 10 symmetric passport $[3^24/3^24/2^34]$
of the tilings \cite[(3.8), (3.11), (3.12)]{Davey:2009bp}
seems to be appropriate for the correspondence construction,
but there are no correspondent passports.
In particular, the genus 0 passport $[3^24/2^41^2/4^22]$
and a presumed quadratic covering branching above the four points $[.../...1^2/...4\,2]$
are not suitable by the reasons discussed in Remark \ref{rm:parity2}.
Accordingly, the tilings \cite[(3.8), (3.11), (3.12)]{Davey:2009bp} 
have no central anti-symmetries.
\end{remark}

\subsection{Maps of degree 12}
\label{sec:deg12b}

At last, here we consider genus 1 Belyi maps with the 
passports  $[3^2\,6/3^2\,6/2^6]$ and $[5\,4\,3/5\,4\,3/2^6]$. 
 The correspondent passports are 
 $[3^2\,6/2^4\,1^4/4^3]$, $[3^2\,6/2^6/2^4\,4]$ 
 and respectively $[5\,4\,3/2^41^4/4^3]$, $[5\,4\,3/2^6/2^4\,4]$. 
There are no genus 0 maps with the last passport.

The genus 0 passport $[3^2\,6/2^4\,1^4/4^3]$ gives two Galois orbits, 
defined over $\IQ$ and $\IQ(\sqrt{-2})$.
Both corresponding genus 1 maps are obtained in \cite{Vidunas:2016xun}
as compositions involving genus 0 Belyi maps.
The $\IQ$-rational map
\begin{align} \label{eq:comb23a} 
\cQ\left\{ \frac{(x-1)^6\,(3x+1)^3\,(x+3)^3}{27\,(x+1)^4\,(x^2+6x+1)^4} \right\} 
\end{align}
is equivalent to the map $\varphi_3$ in \cite[(2.6)]{Vidunas:2016xun} on $y^2=x(x^2-4x+1)$.
One can transform (\ref{eq:comb23a}) to $\varphi_3/(\varphi_3-1)$
by $x\mapsto \big(x^2+2\sqrt{-6}\,y-10x+1\big)/(x+1)^2$.
The dessin is depicted in \cite[(3.32)]{Davey:2009bp} 
and \cite[Fig.~3\refpart{a}]{Vidunas:2016xun}

The $\IQ(\sqrt{-2})$-map
\begin{align}
\cQ\left\{ \frac{256\,x^3\,(x+1)^3\,(x-3)^6}
{\big(4x^3+(21+9\sqrt{-2})x^2+(-6\sqrt{-2}+6)x+\sqrt{-2}+5\big)^4}
 \right\}
\end{align}
is on the elliptic curve
\begin{align}
y^2=16x^3+(63\sqrt{-2}-24)x^2+(-42\sqrt{-2}+96)x+7\sqrt{-2}+8.
\end{align}
The relation to \cite[(5.5)]{Vidunas:2016xun} is via $x\mapsto -3\sqrt{-2}x-\sqrt{-2}+1$.
The dessin is depicted in \cite[(3.35)]{Davey:2009bp} 
and \cite[Appendix \refpart{L}]{Vidunas:2016xun}.


\begin{remark} \label{rm:deg12g00} 
The other correspondent passport $[3^2\,6/2^6/2^4\,4]$
gives  one Belyi map of genus 0:
\begin{align}
\varphi^*_{12}=-\frac{4(x^2+1)^3}{x^4(x^4+3x^2+3)^2}.
\end{align}
However, the map
\begin{align}
\cQ\{\varphi^*_{12}\}=-\frac{(x^2+1)^3}{x^2(x^4+3x^2+3)}
\end{align}
is minimally of genus 0 and degree 6, of the passport $[3^2/6/1^42]$.
As a genus 1 map on $y^2=x^4+3x^2+3$, it has the passport $[3^4/6^2/2^6]$
rather than $[3^2\,6/3^2\,6/2^6]$.
The dessins are in Appendix figure \refpart{O}. 
The genus 1 dessin is 2-isogenous to \cite[Fig.~8\refpart{a}]{Vidunas:2016xun}
\end{remark}

Genus 1 Belyi maps with the passport $[5\,4\,3/5\,4\,3/2^6]$ 
are obtained by the correspondence construction from genus 0 maps
with the passport $[5\,4\,3/2^41^4/4^3]$.
This gives a Galois orbit defined over the degree 6 number field
\begin{equation}
\IQ(\eta)/(\eta^6-2\eta^5-4\eta^4+12\eta^3-14\eta^2+8\eta-4).
\end{equation}
The genus 0 Belyi function is
\begin{align} \label{eq:psid12}
\psi_{12} 
&=
\frac{9x^4\,( 2x+\eta^5-2\eta^4-5\eta^3+16\eta^2-14\eta+4)^5\,F^3}
{200\,(2\eta^5-4\eta^4-5\eta^3+9\eta^2-6\eta+4)\,G^4} \,,
\end{align}
where
\begin{align*} 
F = &\; (\eta^4-4\eta^2-2)\,x+\eta^5+2\eta^4-12\eta^3+2\eta^2-4\eta \,, \\
G = &\; 2x^3+6\eta^2x^2+(-3\eta^5+3\eta^4+18\eta^3-24\eta^2+12\eta)\,x  \\
& +6\eta^5-24\eta^4+36\eta^3-36\eta^2+24\eta-16.
\end{align*}
The function $\cQ\{\psi_{12}\}$ is defined on a genus 1 curve isomorphic to
\begin{align}
y^2= &\;
x^3+3(17\eta^5-58\eta^4-36\eta^3+334\eta^2-462\eta+229)\,x \nonumber \\
& +407\eta^5+43\eta^4-3292\eta^3+890\eta^2+5122\eta-6098.
\end{align}
The dessins are identified as \cite[(3.29)--(3.31),\,(3.33)]{Davey:2009bp}.
They are depicted in appendix figures \refpart{I1}--\refpart{I4}.


\section{Algebraic Painlev\'e VI solutions}\setall
\label{sec:painleve}

The correspondence construction of genus 1 Belyi maps
was inspired by a customary representation of algebraic Painlev\'e VI solutions,
originally due to Kitaev \cite{Kit1} and Boalch \cite{Boa1}. 
The Painlev\'e VI equation is a classical non-linear differential equation 
with the Painlev\'e property: the only movable singularities of its solutions 
are poles \cite{Cla32}.  The recent interest in Painlev\'e equations is mostly 
inspirited by the relation to integrable, isomonodromic differential systems  
\cite{JM81}, \cite{Oka}.
The general form of the  Painlev\'e VI equation is
\begin{align}  \label{eq:pvi}
\frac{d^2q}{dt^2}=\,& {1\over2} \! \left({1\over q}+{1\over q-1}+{1\over q-t}\right) \!
\left(\frac{dq}{dt}\right)^{\!2} \!
-\left({1\over t}+{1\over t-1}+{1\over q-t}\right)\frac{dq}{dt} \nonumber \\
&+{q(q-1)(q-t)\over t^2(t-1)^2}\left(\alpha+\beta {t\over q^2}+
\gamma{t-1\over(q-1)^2}+\delta {t(t-1)\over(q-t)^2}\right) \ .
\end{align} 
When the parameters $\alpha,\beta,\gamma,\delta$ have particular values,
let $\Pvi(\alpha,\beta,\gamma,\delta)$ denote 
the specialized Painlev\'e VI
equation. 

\subsection{An algebraic solution of degree 9}

The algebraic Painlev\'e VI solutions had been recently classified by
Lisovyy and Tykhyy  \cite{LiTy}. Apart from infinite families of rational 
and Picard's $\Pvi(0,0,0,1/2)$ solutions,
there are 45 algebraic solutions up to so-called Okamoto transformations \cite{Oka}.
Most of the interesting cases were originally found  by Boalch \cite{Boa1}, \cite{Boa2}.

In general, it is not practical to present an algebraic Painlev\'e VI function $q(t)$
by its polynomial equation $P(q,t)=0$. 
Let $E$ denote the algebraic curve defined by $P(q,t)=0$. If it 
is of genus 0, a parametrization $(q(s),t(s))$ of $E$ 
is usually a much more compact presentation of the algebraic function $q(t)$.
When the curve $E$ has a positive genus, 
one parametrizes by a simple model of $E$,
thereby considering $q:E\to\IC$ and $t:E\to\IC$ as functions on the underlying curve.
In a minimal parametrization,
the function  $t:E\to\IC$ must be a Belyi map by the Painlev\'e property of $q(t)$.

The solution \#13 in \cite{LiTy} of $\Pvi(1/18,-2/25,2/25,21/50)$ is parametrized as 
\begin{align}
q_{13}=&\; \frac12+\frac{s\,(2s^3+6s^2-63s-350)}{30\,w\,(s+2)},\\
\label{eq:t13}
t_{13}=&\;\frac12+\frac{(s+8)(2s^4-8s^3-42s^2-170s-25)}{54\,w\,(4s+5)},
\end{align}
with $w=\sqrt{s\,(s+8)(4s+5)}$,
thereby defining an elliptic curve whose $j$-invariant equals $889^3/270^2$. 
The map $t_{13}(s,w)$ is indeed a genus 1 Belyi map,
with the passport $[5\,3\,1/5\,3\,1/5\,3\,1]$. 
The particular formula (\ref{eq:t13}) 
matches  the generic $\cQ$-expression (\ref{eq:basicop}).
Indeed $t_{13}=\cQ\{r_{13}\}$,
where
\begin{equation}
r_{13}= 
-\frac{4(s-1)^5(s+5)^3(s-10)}{729\,s\,(4s+5)^3} 
\end{equation}
is a genus 0 map with the passport $[5\,3\,1/2^41/5\,3\,1]$.
The genus 0 and 1 dessins defining the Belyi map $t_{13}(s)$
are depicted in Appendix figure \refpart{J}. The cells are labeled 
by their valency 1, 3 or 5.

Even though $q_{13}(s)$ is not a Belyi map, we can write
\begin{equation}
q_{13}=\cQ\left\{-\frac{4(s-1)^2(s+5)^2(s-10)(s^2+8s+36)}{225(s+8)(4s+5)(s+2)^2}\right\}.
\end{equation}
In recovering the Painleve VI solution $q_{13}(t_{13})$ from these $\cQ$-expressions, 
we must be careful with choosing matching branches of $w$. 
The $\cQ$-forms of $q_{13}$, $t_{13}$ reflect the symmetry 
$\beta+\gamma=0$, $(q,t)\mapsto (1-q,1-t)$ between the singularities $t=0$ and $t=1$. 

The similar solution \#14 in \cite{LiTy} of $\Pvi(2/9,-1/50,1/50,12/25)$ 
is parametrized  with the same $t_{14}=t_{13}$ and 
\begin{equation}
q_{14}=\cQ\left\{-\frac{4(s-1)^4(s+5)^2(s-10)(s+2)}{(4s+5)(s^3+60s+20)^2}\right\}.
\end{equation}
This solution $q_{14}(t_{14})$ was first find by Kitaev \cite{Kit1}.
In total, 17 out of the 45 exceptional solutions in  \cite{LiTy} have genus 1,
while 22 have genus 0.

\subsection{More Belyi maps with the same passport}
\label{sec:notp6}

There are more Belyi maps with the same genus 1 passport
$[5\,3\,1/5\,3\,1/5\,3\,1]$.
Several of them can be obtained by 
genus 0 maps with the correspondent passport $[5\,3\,1/2^41/5\,3\,1]$. 
First of all, the dual rational map $1/r_{13}$ leads to a non-equivalent genus 1 map:
\begin{align} \label{eq:t13x}
t^*_{13}= &\; \cQ\left\{\frac1{r_{13}}\right\} \\
= &\; \frac12+\frac{(s+8)(2s^4-8s^3-42s^2-170s-25)}
{8\,(s-1)^2\,(s+5)\,\sqrt{(s-1)(s+5)(s+8)(s-10)}}. \nonumber
\end{align}
The $j$-invariant of the elliptic curve defined by the radical, viz.,
$w^2 = (s-1)(s+5)(s+8)(s-10)$, equals $2\cdot 42^3/5^2$.
The genus 0 and 1 dessins are depicted in Appendix figure \refpart{K}.


Further, there is a Galois orbit of degree 5. The field of definition is
\begin{equation}
\IQ(\zeta)/(\zeta^5-2\zeta^3-4\zeta^2-6\zeta-4)
\;. 
\end{equation}
The genus 1 Belyi map is
\begin{align} \label{eq:p6d5}
\cQ\left\{ \frac{(x+1)^5\,(5x-c_1)^3\,(10x-c_2)}{(x-1)^5\,(5x+c_1)^3\,(10x+c_2)}
\right\}
\end{align}
with $c_1=13\zeta^4-6\zeta^3-46\zeta^2-24\zeta+13$, $c_2=3\zeta^4-\zeta^3-6\zeta^2-9\zeta-12$.
The dessin orbit is given in \refpart{L1}--\refpart{L3}. There are no other Belyi maps
with the passport $[5\,3\,1/2^41/5\,3\,1]$, but there might be genus 1 maps
with the passport $[5\,3\,1/5\,3\,1/5\,3\,1]$ with no symmetries between 
the three fibers. 

Dessins 
exists that are not obtained by the correspondence construction.
For example, we can modify the Appendix tiling \refpart{L2} topologically
by flipping a loose edge to another cell without affecting the cell valencies,
so that a tiling on the right-hand side 
of Appendix figure 
\refpart{M} is obtained. that tiling has the same passport, but no central anti-symmetry.

\appendix
\section{Appendix: Explicit Dessins}\setall

In this appendix, we collect the various dessins which have been discussed in the text.
We present (1) the genus zero dessin on the left; (2) the genus 1 dessin obtained by quadratic correspondence on the right as a tiling of the doubly-periodic plane; (3) the j-invariant of the elliptic curve for the genus 1 dessin; (4) the associated Belyi map, cross-referencing the main text, and, where applicable, (5) the associated brane-tiling in the catalogue of \cite{Davey:2009bp}.
The explicit cuts on the genus 0 dessin are marked in blue.
This also demonstrates a method of getting genus 1 tilings without necessarily computing the Belyi maps.

\noindent
\begin{picture}(440,158)
\put(0,0){\includegraphics[width=440pt]{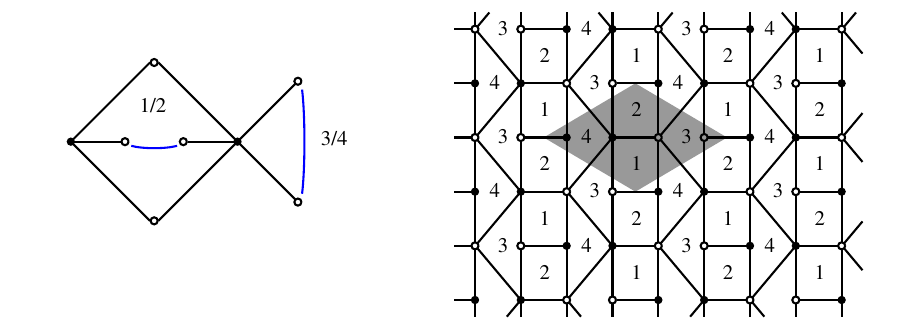}}
\put(0,10){\refpart{A} Map \eqref{eq:f26}, $\!j\!=\!-5000$;
\cite[\!(2.6)]{Davey:2009bp}.}
\end{picture}
\begin{picture}(440,158)
\put(0,0){\includegraphics[width=440pt]{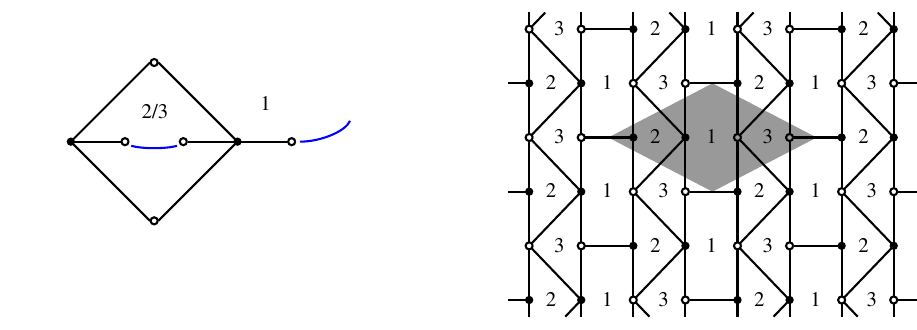}}
\put(0,10){\refpart{B1} Map \eqref{eq:f22},\! with {\small $\sqrt7$};
\cite[\!(2.2)]{Davey:2009bp}, {\em \!SPP}.}
\end{picture}
\begin{picture}(440,158)
\put(0,0){\includegraphics[width=440pt]{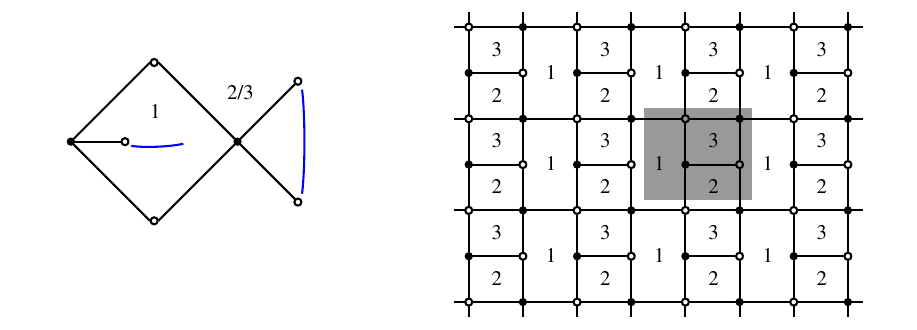}}
\put(0,10){\refpart{B2} {\small $\sqrt{7\,}$}-conjugated \eqref{eq:f22};
\cite[\!(2.3)]{Davey:2009bp}.}
\end{picture}
\begin{picture}(440,158)
\put(0,0){\includegraphics[width=440pt]{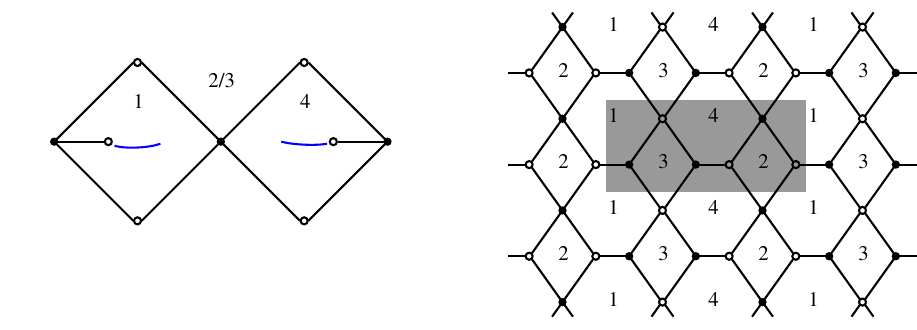}}
\put(0,10){\refpart{C} Map \eqref{eq:comb33a}; 
\cite[\!(3.5)]{Davey:2009bp}, $\!L^{222}$(II).}
\end{picture}
\begin{picture}(440,158)
\put(0,0){\includegraphics[width=440pt]{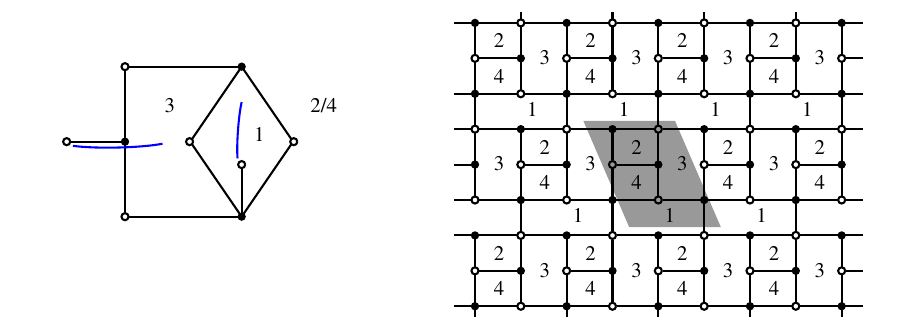}}
\put(0,10){\refpart{D} Map \eqref{eq:qm310},\! 
with {\small$\!\sqrt{-5}$};\! \cite[\!(3.10)]{Davey:2009bp}.}
\end{picture}
\begin{picture}(440,158)
\put(0,0){\includegraphics[width=440pt]{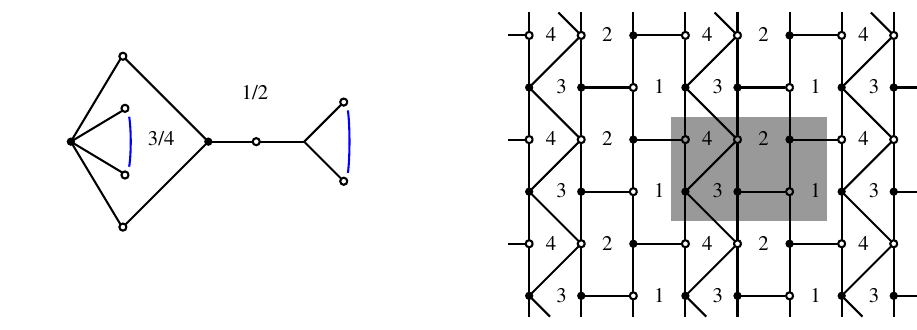}}
\put(0,10){\refpart{E1} Map \eqref{eq:deg10d5}, $\xi\!\in\!\IR;\!$
\cite[\!(3.4)]{Davey:2009bp}, $\!L^{131}$.}
\end{picture}
\begin{picture}(440,158)
\put(0,0){\includegraphics[width=440pt]{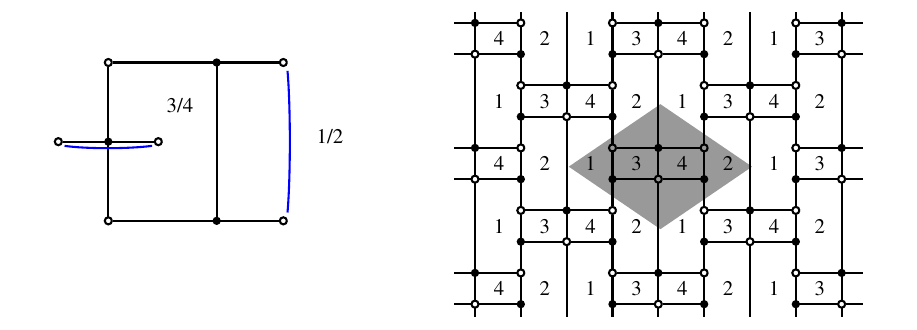}}
\put(0,10){\refpart{E2} Map \!\eqref{eq:deg10d5}, $\xi\!\in\!\IR;\!$
\cite[\!(3.6)]{Davey:2009bp}, $\!dP_1$.}
\end{picture}
\begin{picture}(440,158)
\put(0,0){\includegraphics[width=440pt]{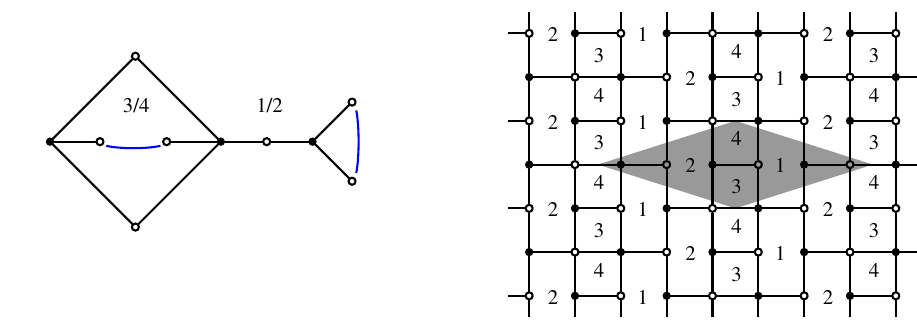}}
\put(0,10){\refpart{E3} Map \eqref{eq:deg10d5}, $\xi\in\IR$;
\cite[\!(3.7)]{Davey:2009bp}.}
\end{picture}
\begin{picture}(440,158)
\put(0,0){\includegraphics[width=440pt]{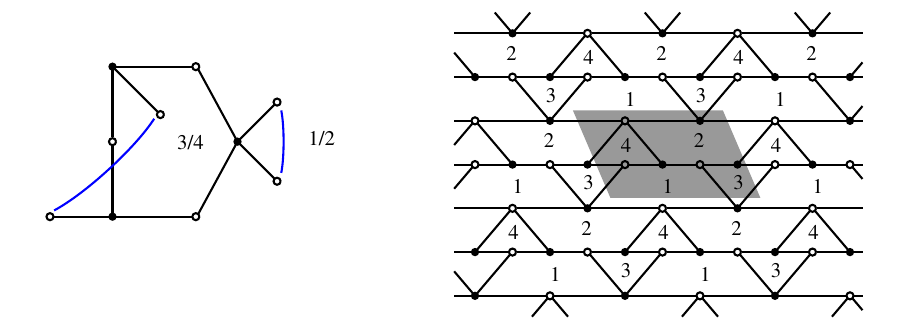}}
\put(0,10){\refpart{E4} Map \eqref{eq:deg10d5}, $\xi\not \in\IR$;
\cite[\!(3.9)]{Davey:2009bp}.}
\end{picture}
\begin{picture}(440,158)
\put(0,0){\includegraphics[width=440pt]{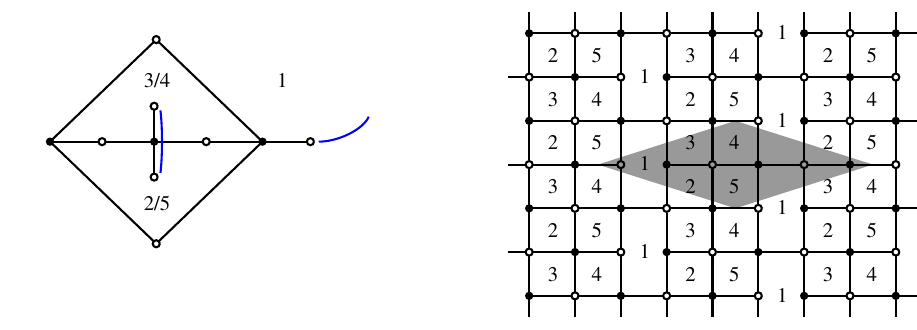}}
\put(0,10){\refpart{F1} Map \eqref{eq:deg11r11}; 
\cite[\!(3.19)]{Davey:2009bp}.}
\end{picture}
\begin{picture}(440,158)
\put(0,0){\includegraphics[width=440pt]{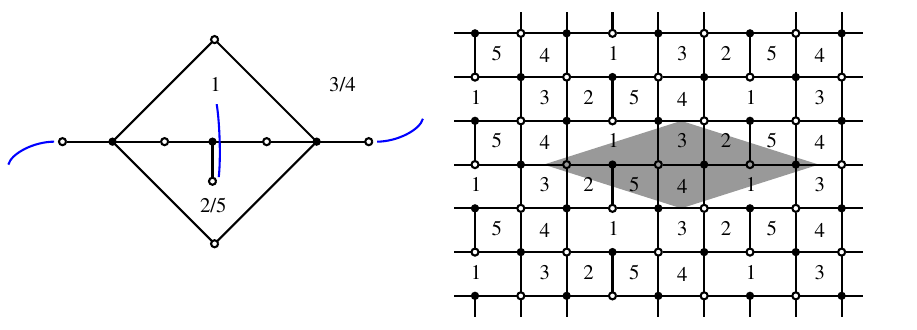}}
\put(0,10){\refpart{F2} \small{$\!\sqrt{11\,}$}-conjugated \eqref{eq:deg11r11};\!
\cite[\!(3.18)]{Davey:2009bp}.}
\end{picture}
\begin{picture}(440,158)
\put(0,0){\includegraphics[width=440pt]{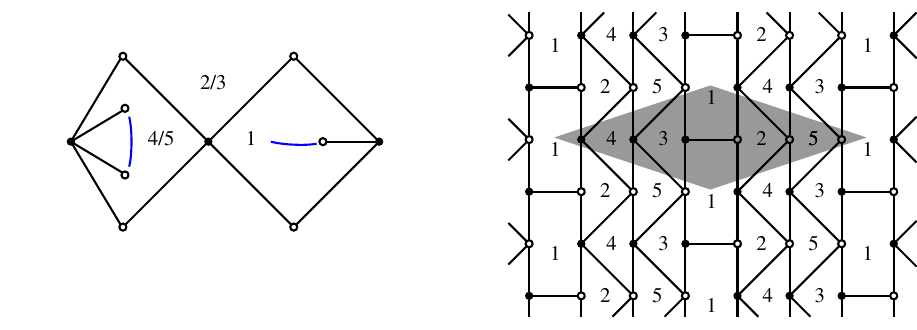}}
\put(0,10){\refpart{G1} Map \eqref{eq:deg11cubic}, $\mu\!\in\!\IR$;
\cite[\!(3.13)]{Davey:2009bp}, $\!L^{232}$.}
\end{picture}
\begin{picture}(440,158)
\put(0,0){\includegraphics[width=440pt]{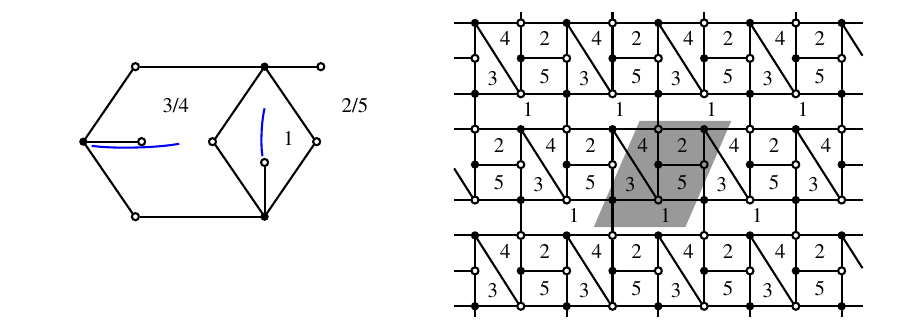}}
\put(0,10){\refpart{G2} Map \eqref{eq:deg11cubic}, $\mu\not\in\IR$;
\cite[\!(3.20)]{Davey:2009bp}.}
\end{picture}
\begin{picture}(440,158)
\put(0,0){\includegraphics[width=440pt]{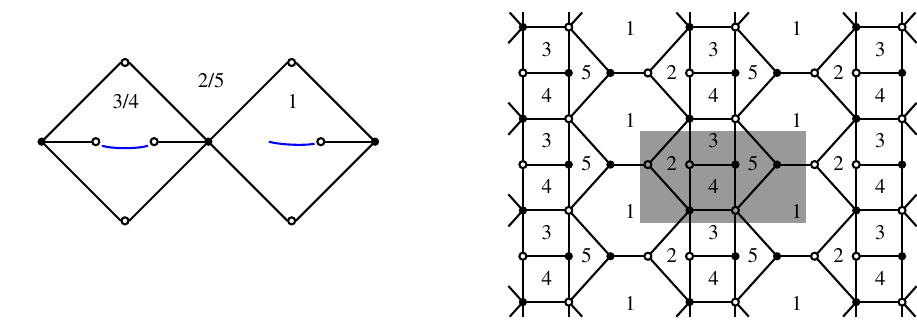}}
\put(0,10){\refpart{H1} Map \eqref{eq:deg10d5}, $\lambda\in\IR$;
\cite[\!(3.16)]{Davey:2009bp}.}
\end{picture}
\begin{picture}(440,158)
\put(0,0){\includegraphics[width=440pt]{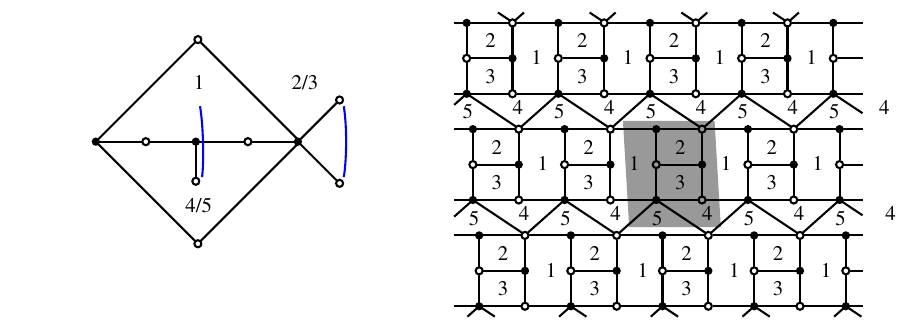}}
\put(0,10){\refpart{H2} Map \eqref{eq:deg10d5}, $\lambda\not\in\IR$;
\cite[\!(3.21)]{Davey:2009bp}.}
\end{picture}
\begin{picture}(440,158)
\put(0,0){\includegraphics[width=440pt]{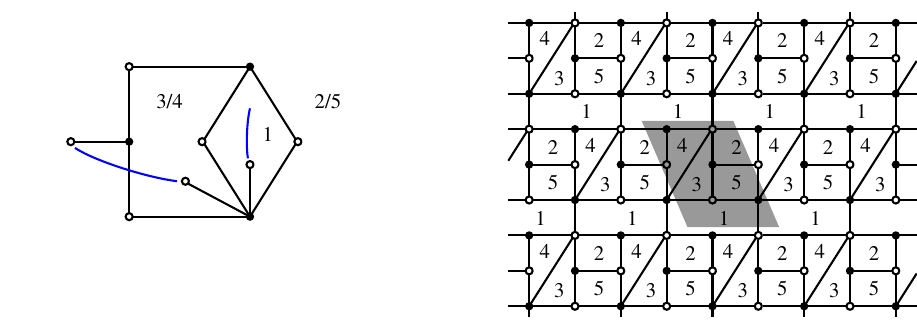}}
\put(0,10){\refpart{H3} Map \eqref{eq:deg10d5}, $\lambda\not\in\IR$;
\cite[\!(3.24)]{Davey:2009bp}.}
\end{picture}
\begin{picture}(440,158)
\put(0,0){\includegraphics[width=440pt]{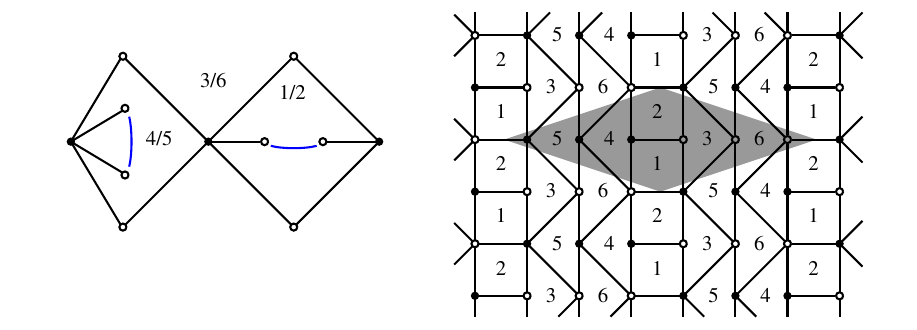}}
\put(0,10){\refpart{I1} Map \eqref{eq:psid12}, $\eta\in\IR$;
\cite[\!(3.30)]{Davey:2009bp}.}
\end{picture}
\begin{picture}(440,158)
\put(0,0){\includegraphics[width=440pt]{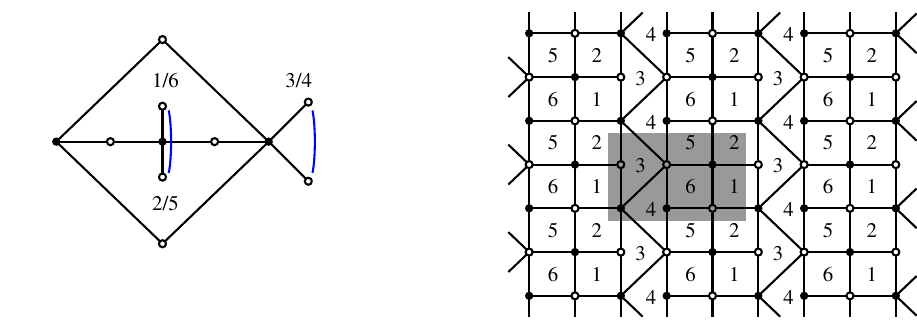}}
\put(0,10){\refpart{I2} Map \eqref{eq:psid12}, $\eta\in\IR$;
\cite[\!(3.31)]{Davey:2009bp}.}
\end{picture}
\begin{picture}(440,158)
\put(0,0){\includegraphics[width=440pt]{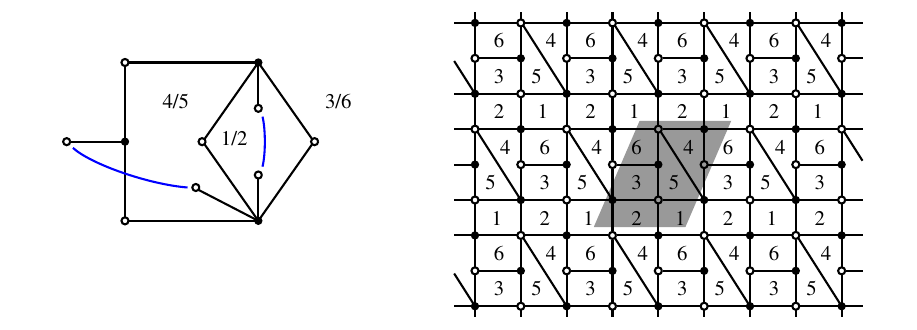}}
\put(0,10){\refpart{I3} Map \eqref{eq:psid12}, $\eta\not\in\IR$;
\cite[\!(3.29)]{Davey:2009bp}.}
\end{picture}
\begin{picture}(440,158)
\put(0,0){\includegraphics[width=440pt]{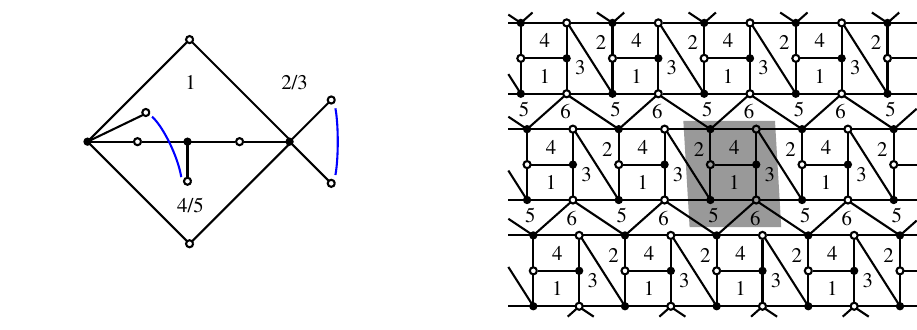}}
\put(0,10){\refpart{I4} Map \eqref{eq:psid12},  $\eta\not\in\IR$;
\cite[\!(3.33)]{Davey:2009bp}.}
\end{picture}
\begin{picture}(440,158)
\put(0,0){\includegraphics[width=440pt]{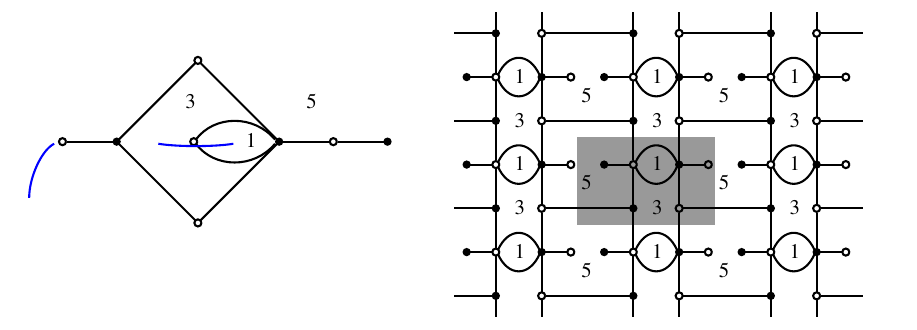}}
\put(0,10){\refpart{J} Dessin for Painlev\'e VI solution\! (\ref{eq:t13}).}
\end{picture}
\begin{picture}(440,158)
\put(0,0){\includegraphics[width=440pt]{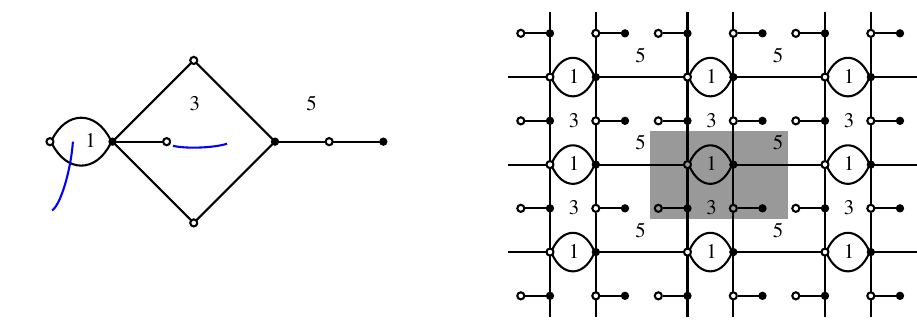}}
\put(0,10){\refpart{K} Dessin for function (\ref{eq:t13x}).}
\end{picture}
\begin{picture}(440,158)
\put(0,0){\includegraphics[width=440pt]{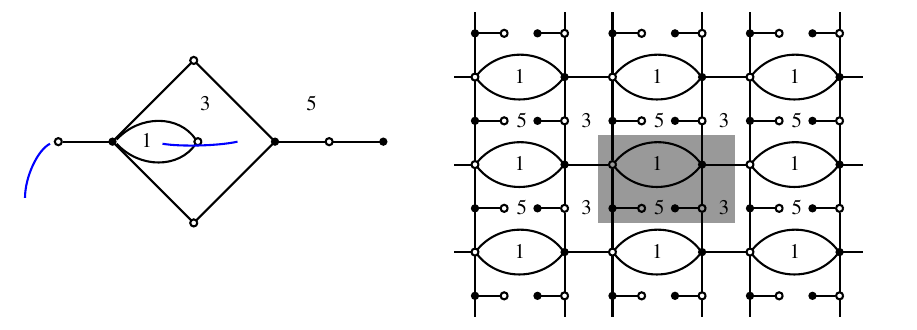}}
\put(0,10){\refpart{L1} Dessin for (\ref{eq:p6d5}), with $\zeta\in\IR$.}
\end{picture}
\begin{picture}(440,158)
\put(0,0){\includegraphics[width=440pt]{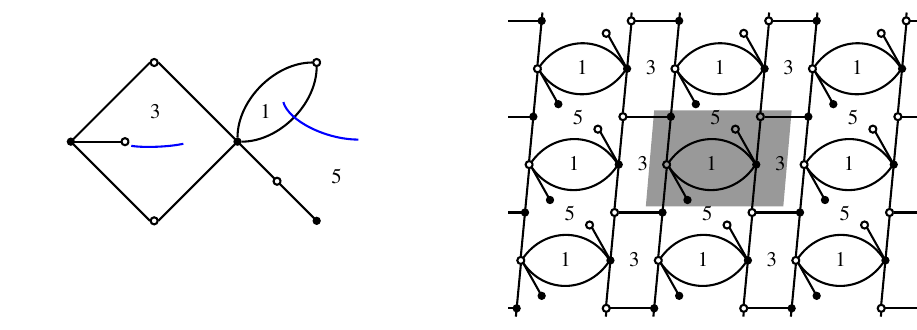}}
\put(0,10){\refpart{L2} Dessin for (\ref{eq:p6d5}), with $\zeta\not\in\IR$.}
\end{picture}
\begin{picture}(440,158)
\put(0,0){\includegraphics[width=440pt]{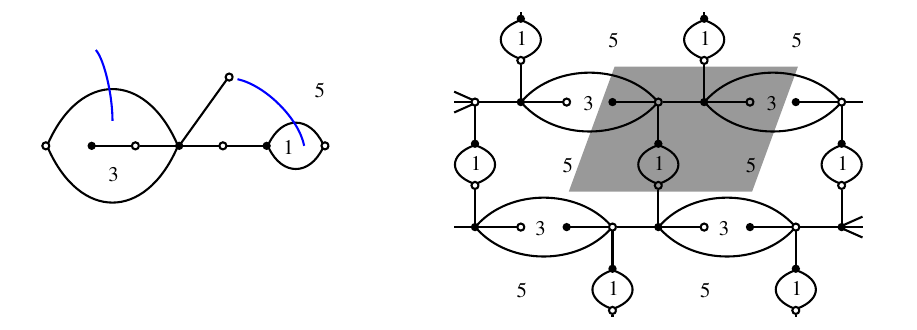}}
\put(0,10){\refpart{L3} Dessin for (\ref{eq:p6d5}), with $\zeta\not\in\IR$.}
\end{picture}
\begin{picture}(440,184)
\put(0,26){\includegraphics[width=440pt]{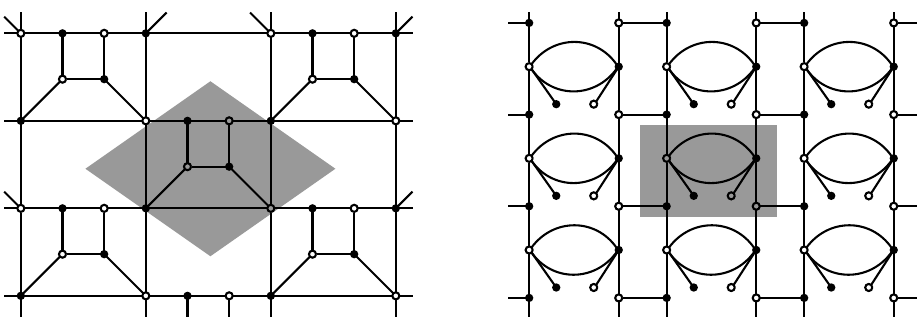}}
\put(0,10){\!\refpart{M} Dessins not obtainable by the quadratic correspondence;\! 
\cite[\!(3.25)]{Davey:2009bp}\! on the left.}
\end{picture}
\begin{picture}(440,158)
\put(0,0){\includegraphics[width=440pt]{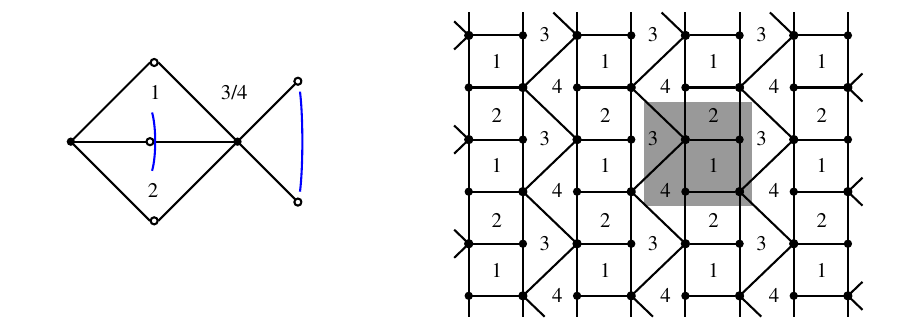}}
\put(0,10){\refpart{N} Non-bipartite tiling in Remark \ref{rm:parity2}.}
\end{picture}
\begin{picture}(440,158)
\put(0,0){\includegraphics[width=440pt]{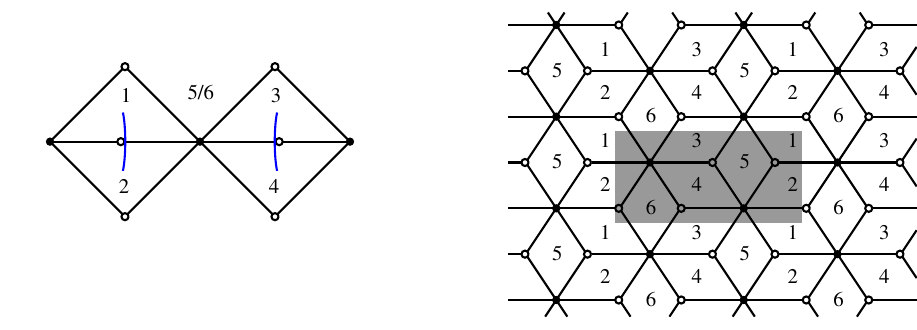}}
\put(0,10){\refpart{O} Example for Remark \ref{rm:deg12g00}.}
\end{picture}

\vspace{1in}
\section*{Acknowledgments}
We are most grateful to the organizers for inviting us to the {\it Workshop on Grothendieck-Teichm\"uller Theories}, 2016, at the Chern Institute of Mathematics, Nankai University where the most congenial and productive atmosphere inspired the beginnings of our friendship and programme.

YHH would like to
thank the Science and Technology Facilities Council, UK, for grant ST/J00037X/1, the Chinese
Ministry of Education, for a Chang-Jiang Chair Professorship at NanKai University as well as
the City of Tian-Jin for a Qian-Ren Scholarship, and Merton College, Oxford, for her enduring
support.

\small

\end{document}